\documentclass{article}
\usepackage{amsmath}
\usepackage{amssymb}
\usepackage{amsthm}
\usepackage{enumerate}
\usepackage{mathrsfs}
\usepackage[multiple]{footmisc}
\usepackage[a4paper, total={6in, 9in}]{geometry}

\title{The star avoidance game}
\author{Adrian Beker \footnote{Trinity College, Cambridge CB21TQ, UK} \footnote{\href{mailto:ab2454@cam.ac.uk}{\nolinkurl{ab2454@cam.ac.uk}}}}
\date{\today}

\begin{document}

\maketitle

\begin{abstract}
    Let $n, k$ be positive integers. The $(k+1)$-star avoidance game on $K_n$ is played as follows. Two players take it in turn to claim a (previously unclaimed) edge of the complete graph on $n$ vertices. The first player to claim all edges of a subgraph isomorphic to a $(k+1)$-star loses. Equivalently, each player must keep all degrees in the subgraph formed by his edges at most $k$. If all edges have been chosen and neither player has lost, the game is declared a draw. We prove that, for each fixed $k$, the game is a win for the second player for all $n$ sufficiently large.
\end{abstract}

\section{Introduction}

Many natural combinatorial games occur as follows. We have a finite set (called the \emph{board}), some subsets of which are designated as \emph{lines}. Two players take it in turn to claim a (previously unclaimed) element of the board. The first player to complete a line loses (and the other player is declared the winner). If all elements have been chosen and neither player has lost, the game is declared a draw. Due to the winning criterion, games of the described kind are called \emph{mis\`ere} games. Games with the usual winning criterion of making the desired object are called `achievement games' -- see Beck \hyperlink{beck1}{[1]} for a discussion of both kinds of game and Slany \hyperlink{slany}{[4]} for background on mis\`ere games. For a related more general overview of results and methods in combinatorial game theory, see Beck \hyperlink{beck2}{[2]}.

A subfamily of mis\`ere games of particular interest is the class of \emph{sim-like} games. These games have board the edge set of the complete graph $K_n$. The lines are subsets which form a subgraph isomorphic to a graph from some fixed family $\mathcal{F}$. The first example of such a game we are aware of is the game of \emph{Sim}, which is played on $K_6$ and has $\mathcal{F} = \{K_3\}$ (see Simmons \hyperlink{simmons}{[3]}). For a survey of sim-like games from a more computational perspective, see Slany \hyperlink{slany}{[4]}.

It is an easy consequence of Ramsey's theorem that any sim-like game is not a draw for all sufficiently large boards. Indeed, take $G \in \mathcal{F}$ and consider $n \geq R(k, k)$, where $k$ is the order of $G$. Then if all edges of $K_n$ have been chosen, at least one player has claimed all edges of a $k$-clique and hence of a subgraph isomorphic to $G$.

In this paper we will consider the sim-like game given by $\mathcal{F} = \{S_{k + 1}\}$, where $k$ is a fixed positive integer (and $S_{k + 1}$ denotes the graph $K_{1, k + 1}$), also called the $(k + 1)$-\emph{star avoidance game}. This is one of the simplest and most natural mis\`ere games on a graph. Our main result is the following:
\\\\
\textbf{Theorem 1.} \textit{The second player wins the $(k+1)$-star avoidance game on $K_n$ whenever $n \geq 200k$.}\\

In Section 2, we give a proof of Theorem 1. In Section 3, we conclude the paper with some remarks and a discussion of related open problems.
Throughout the paper, all graphs in consideration are simple and we assume standard notation from graph theory. We abbreviate the first and second player to PI and PII respectively. A 'round' comprises a move of PI followed by a move of PII.

\newpage

\section{Proof of Theorem 1}

The aim of this section is to provide a proof of the following theorem:
\\\\
\textbf{Theorem 1.} \textit{PII wins the $(k+1)$-star avoidance game on $K_n$ whenever $n \geq 200k$.}\\

We first give a brief overview of our approach. In the context of the $(k+1)$-star avoidance game, by 'valid subgraph' we will mean 'subgraph of $K_n$ of maximum degree at most $k$'. We note that a straightforward way for PII to win would be to build a valid subgraph of size ex$(n, S_{k+1}) = \left\lfloor \frac{nk}{2} \right\rfloor$. However, in the case when $nk$ is even, this approach would require careful adjustments to the opponent's actions in the final stage of the game, and it turns out that in fact slightly less is needed. Instead, we contend that PII can build a subgraph of size $\left\lfloor \frac{nk - 1}{2} \right\rfloor$ with the property that, in the case when $nk$ is even, there exists an unclaimed edge that extends it to a valid subgraph. We then use the fact that PI's last move is uniquely determined to argue that PII can force a win.

In order to prove the main claim, we define an auxiliary game, called the \emph{pair clipping game}, which is strictly speaking not a positional game. However, the strategy required for PII to win can be viewed as very similar to a fast winning strategy in the Maker-Breaker perfect matching game. The difference is in the winning criterion, which is slightly modified according to the needs of our claim. Furthermore, we require this game to be played on a general graph instead of $K_n$. We make use of a winning strategy in the pair clipping game by building 'layers' of almost perfect matchings, until we reach PII's goal.

Finally, the main part of the proof deals with finding a strategy for the pair clipping game, which is provided by Theorem 5. Note that, by the way the game is defined, the graph represents 'obstacles' for PII, i.e. previously claimed edges in the context of the star avoidance game. PII essentially employs a greedy strategy, in which he tries to keep the graph as sparse as possible. It entails inductively controlling both the average and maximum degree of the graph.\\

We start by defining the auxiliary game. Let $G$ be a graph and let $n = v(G)$. The pair clipping game on $G$, PCG$(G)$, is defined as follows. It is played by two players and consists of rounds of the following form:
\begin{itemize}
    \item PI adds at most one edge which is not already present in $G$ and then PII removes two non-adjacent vertices from $G$.
\end{itemize}

PII wins the game if he can make $\lfloor \frac{n - 1}{2} \rfloor$ moves so that $G$ becomes empty. For convenience, we will denote by $G_j$ the graph $G$ after $j$ turns. Thus, $G_0$ is the initial graph and $v(G_{2j}) = v(G_{2j+1}) = v(G_0) - 2j$ holds for all $j \geq 0$. Likewise, we will denote by $\{u_j, v_j\}$ the pair of vertices chosen in the $j$-th turn (or an arbitrary element of $V(G_{j-1})^{(2)}$ if PI does nothing in the $j$-th turn).\\

We define the notion of a nice pair of vertices, which is central to the greedy strategy. Let $G$ be a graph and let $u, v \in V(G)$ be distinct. We say that the pair $\{u, v\}$ is \emph{nice} if $\{u, v\} \not \in E(G)$ and $d_G(u) + d_G(v) \geq 2d(G)$, i.e. $\{u, v\}$ is not an edge and the average of the degrees of $u$ and $v$ is at least the average degree of the whole graph. The following lemma guarantees the existence of such a pair in graphs with not too large maximum degree:
\\\\
\textbf{Lemma 2.} \textit{Let $G$ be a graph with $v(G) \geq 2$ and $\Delta(G) \leq \frac{1}{2}v(G) - 1$. Then $G$ has a nice pair.}
\\\\
\textit{Proof.} If $v(G) = 2$, then we are done since $G$ is empty, so suppose $v(G) \geq 3$. Let $H$ be the complement of $G$ and let $n = v(G) = v(H)$. Then note that $\delta(H) = n - 1 - \Delta(G) \geq \frac{1}{2}n$, so by Dirac's theorem, $H$ has a Hamiltonian cycle $v_1v_2\ldots v_n$. Averaging over this cycle, we obtain $$\frac{1}{n}\sum_{j = 1}^{n} \frac{d_H(v_j) + d_H(v_{j+1})}{2} = \frac{1}{n}\sum_{j = 1}^{n} d_H(v_j) = d(H).$$ Hence, there exists $j \in [n]$ such that $\frac{d_H(v_j) + d_H(v_{j+1})}{2} \leq d(H)$. But note that $\{v_j, v_{j + 1}\} \not \in E(G)$ and $$d_G(v_j) + d_G(v_{j+1}) = (n - 1 - d_H(v_j)) + (n - 1 - d_H(v_{j+1})) \geq 2(n - 1) - 2d(H) = 2d(G),$$ so $\{v_j, v_{j+1}\}$ is a nice pair in $G$. $\blacksquare$\\

We now introduce a certain notion of sparseness of graphs which will be used in Theorem 5. Let $f, g : \mathbb{N} \to \mathbb{R}_{\geq 0}$ be functions. Given a graph $G$, we say $G$ is \emph{$g$-sparse} if $d(G) \leq g(v(G))$. If $G$ additionally satisfies $\Delta(G) \leq f(v(G))$, we say $G$ is \emph{$(f, g)$-sparse}. The following lemma is mostly technical and describes the effect of removing a nice pair on the sparseness of a graph: 
\\\\
\textbf{Lemma 3.} \textit{Let $\alpha > 0$ and let $g : \mathbb{N} \to \mathbb{R}_{\geq 0}$ be given by $g(n) = \alpha n + 1$. Let $G$ be a $g$-sparse graph with $v(G) \geq 4$ and suppose that $G'$ is obtained from $G$ by adding at most one edge. Let $G'' = G'[V(G') \setminus \{u, v\}]$, where $\{u, v\} \not \in E(G')$. If $d_{G'}(u) + d_{G'}(v) \geq 2d(G)$, then $G''$ is $g$-sparse. In particular, if $\{u, v\}$ is nice in $G'$, then $G''$ is $g$-sparse.}
\\\\
\textit{Proof.} Observe that 
\begin{equation*} 
\begin{split}
d(G'') &= \frac{2e(G'')}{v(G'')} = \frac{2e(G') - 2(d_{G'}(u) + d_{G'}(v))}{v(G') - 2} \leq \frac{v(G)d(G) + 2 - 2 \cdot 2d(G)}{v(G) - 2}\\
&= \frac{(v(G) - 4)d(G) + 2}{v(G) - 2} \leq \frac{(v(G) - 4)(\alpha v(G) + 1) + 2}{v(G) - 2} = \frac{\alpha v(G)(v(G) - 4)}{v(G) - 2} + 1\\
&< \alpha(v(G) - 2) + 1 = \alpha v(G'') + 1,
\end{split}
\end{equation*}
so $G''$ is $g$-sparse. If $\{u, v\}$ is nice in $G'$, we are done since $d(G') \geq d(G)$. $\blacksquare$ 
\\\\
The next lemma is straightforward and serves mainly for the base cases in the proof of Theorem 5:
\\\\
\textbf{Lemma 4.} \textit{Let $G$ be a graph. Then PII wins PCG($G$) if either
\begin{enumerate}[(i)]
    \item $v(G) \geq 3$ and $G$ is $(1, 1)$-sparse or
    \item $v(G) \geq 5$ and $G$ is $1$-sparse.
\end{enumerate}}
\noindent \textit{Proof.} To prove part $(i)$, we use induction on the order of $G$. If $v(G) \in \{3, 4\}$, then $\Delta(G) \leq 1$, so it is easy to see that PII wins. If $v(G) \geq 5$, then PII removes $u_1$ and any non-adjacent vertex -- he can do so since $d_{G_1}(u_1) \leq 2 \leq v(G_1) - 2$. In this way, $G_2$ is $(1, 1)$-sparse, so we are done by the induction hypothesis.

For part $(ii)$, we also use induction on $v(G)$. Note that we may assume that $\Delta(G) \geq 2$ since otherwise $G$ is $(1, 1)$-sparse and we are done by part $(i)$. If $v(G) \in \{5, 6\}$, then it is easy to check that PII can ensure that $G_2$ is $(1, 1)$-sparse, so part $(i)$ again applies. If $v(G) \geq 7$, then PII removes a vertex of degree at least $2$ in $G_1$ and any non-adjacent vertex -- he can do so since $\Delta(G_1) \leq e(G_1) \leq \frac{v(G_1)}{2} + 1 \leq v(G_1) - 2$. In this way, we have $e(G_2) \leq e(G_1) - 2 \leq \frac{v(G)}{2} - 1 = \frac{v(G_2)}{2}$, so $G_2$ is $1$-sparse, as desired. $\blacksquare$
\\\\
\textbf{Theorem 5.} \textit{Let $f, g : \mathbb{N} \to \mathbb{R}_{\geq 0}$ be given by $f(n) = \frac{n - 1}{2}$, $g(n) = \frac{n}{100} + 1$. Then PII wins PCG$(G)$ for any $(f, g)$-sparse graph $G$.}
\\\\
\textit{Proof.} We proceed by induction on the order of $G$. Suppose that $G$ is a graph with $n = v(G)$ and  $\Delta(G) \leq f(n)$, $d(G) \leq g(n)$. If $n \in \{1, 2\}$, then $G$ is empty, so PII immediately wins. If $n \in \{3, 4\}$, then $G$ is $(1, 1)$-sparse, so we are done by part $(i)$ of Lemma 4. Hence, we may assume that $n \geq 5$. Consider first the case when $n < 10$. Then we have $$v(G)d(G) \leq ng(n) \leq n\left(\frac{n}{100} + 1\right) = n + \frac{n^2}{100} < n + 1.$$ Since $v(G)d(G) = 2e(G)$ is an integer, we in fact have $v(G)d(G) \leq n$, whence $d(G) \leq 1$. Therefore, we are done by part $(ii)$ of Lemma 4. From now on, we assume that $n \geq 10$. We consider two cases:   
\\\\
\textbf{Case 1.} $\Delta(G) \leq f(n - 2) - 1$
\\\\
We have $\Delta(G_1) \leq \Delta(G_0) + 1 \leq f(n - 2) < \frac{1}{2}n - 1$. By Lemma 2, $G_1$ has a nice pair $\{u_2, v_2\}$, so PII removes $\{u_2, v_2\}$. By Lemma 3, $G_2$ is $g$-sparse, and hence $(f, g)$-sparse because $\Delta(G_2) \leq \Delta(G_1)$. Thus, we are done by the induction hypothesis.
\\\\
\textbf{Case 2.} $f(n - 2) - 1 < \Delta(G) \leq f(n)$
\\\\
In the $j$-th round, PII acts as follows:
\begin{itemize}
    \item Remove a vertex of maximum degree in $G_{2j - 1}$ and any non-adjacent vertex.
\end{itemize}
We first make the following easy observations:
\\\\
\textbf{Claim A.} \textit{For all $j \geq 0$ such that PII can make a move in each of the first $j + 1$ rounds, the following hold:
\begin{enumerate}[(i)]
    \item $\Delta(G_{2j}) \leq \Delta(G_{2j+1}) \leq \Delta(G_{2j}) + 1$,
    \item $e(G_{2j+2}) \leq e(G_{2j+1}) - \Delta(G_{2j+1}) \leq e(G_{2j}) - \Delta(G_{2j}) + 1$,
    \item $\Delta(G_{2j+2}) \leq \Delta(G_{2j}).$
\end{enumerate}}

\noindent \textit{Proof.} We note that $(i)$ and $(ii)$ are clear. To see that $(iii)$ holds, note that we are done if $\Delta(G_{2j+1}) = \Delta(G_{2j})$ because $\Delta(G_{2j+2}) \leq \Delta(G_{2j+1})$. However, if $\Delta(G_{2j+1}) = \Delta(G_{2j}) + 1$, then any vertex of maximum degree in $G_{2j+1}$ must be incident to the edge $\{u_{2j+1}, v_{2j+1}\}$. Hence, the conclusion follows. $\blacksquare$
\\\\
Let $r(n) = \left\lfloor \frac{n-1}{4} \right\rfloor$. For all $1 \leq j \leq r(n)$, note that we have $$v(G_{2j-1}) - 2 = n - 2j \geq n - 2r(n) \geq f(n) + 1 \geq \Delta(G_0) + 1 \geq \Delta(G_{2j-2}) + 1 \geq \Delta(G_{2j-1}),$$ so PII can make a move in the $j$-th round. Therefore, by the induction hypothesis, it suffices to show that there exists $j \in [r(n)]$ such that $G_{2j}$ is $(f, g)$-sparse. So suppose for contradiction that this is not the case.
\\\\
\textbf{Claim B.} \textit{For all $1 \leq j \leq r(n)$, the following hold:
\begin{enumerate}[(i)]
    \item $\Delta(G_{2j}) > f(n - 2j)$,
    \item $d(G_{2j}) \leq g(n - 2j)$, i.e. $G_{2j}$ is $g$-sparse.
\end{enumerate}}
\noindent \textit{Proof.} We proceed by induction on $j$. Note that $(i)$ follows from $(ii)$ combined with the assumption that $G_{2j}$ is not $(f, g)$-sparse, so it suffices to prove $(ii)$. Letting $s = d_{G_{2j-1}}(u_{2j}) + d_{G_{2j-1}}(v_{2j})$, we know that $s \geq \Delta(G_{2j-1}) \geq \Delta(G_{2j-2})$. By Lemma 3, it suffices to show that $s \geq 2d(G_{2j - 2})$. To this end, we note that if $j = 1$, then $$s \geq \Delta(G_0) \geq f(n - 2) - \frac{1}{2} \stackrel{(\dagger)}{\geq} 2g(n) \geq 2d(G_0),$$ as desired. Note that the inequality $(\dagger)$ is equivalent to $n \geq \frac{25}{3}$, which indeed holds by assumption. Similarly, if $j > 1$, then we have $$s \geq \Delta(G_{2j-2}) \stackrel{(i)}{\geq} f(n - 2j + 2) + \frac{1}{2} \stackrel{(*)}{\geq} 2g(n - 2j + 2) \stackrel{(ii)}{\geq} 2d(G_{2j-2}),$$ as desired. Note that we used the induction hypothesis in the inequalities $(i)$ and $(ii)$. Moreover, the inequality $(*)$ is equivalent to $n - 2j + 2 \geq \frac{25}{6}$, which holds as $n - 2j + 2 \geq n - 2 \cdot \frac{n - 1}{4} + 2 = \frac{n + 5}{2} \geq \frac{15}{2}$. $\blacksquare$
\\\\ 
Using part $(ii)$ of Claim A and part $(i)$ of Claim B, we obtain
\begin{equation*} 
\begin{split}
e(G_0) &\geq e(G_0) - e(G_{2r(n)}) + \Delta(G_{2r(n)}) = \Delta(G_{2r(n)}) + \sum_{j = 1}^{r(n)} e(G_{2j-2}) - 
e(G_{2j})\\
&\geq \Delta(G_{2r(n)}) + \sum_{j = 1}^{r(n)} (\Delta(G_{2j-2}) - 1) \geq f(n - 2r(n)) + f(n - 2) - 1 + \sum_{j = 2}^{r(n)} \left(f(n - 2j + 2) - \frac{1}{2}\right)\\
&\geq \sum_{j = 1}^{r(n) + 1} \left(\frac{n}{2} - j\right) - 1 = \frac{(r(n) + 1)(n - r(n) - 2)}{2} - 1 \geq \frac{\frac{n}{4}\left(n - \frac{n - 1}{4} - 2\right)}{2} - 1 = \frac{n(3n - 7)}{32} - 1.
\end{split}
\end{equation*}
On the other hand, the assumption on the $g$-sparseness of $G_0$ implies that $$e(G_0) = \frac{v(G_0)d(G_0)}{2} \leq \frac{ng(n)}{2} = \frac{n\left(\frac{n}{100} + 1\right)}{2}.$$ Since $\frac{n(3n - 7)}{32} - 1 > \frac{n\left(\frac{n}{100} + 1\right)}{2}$, we obtain the desired contradiction. $\blacksquare$\\

The following notation will be found useful in the proof of Theorem 1, and in fact applies to any sim-like game. For $j \in \{1, 2\}$, let $H_{j, t}$ be the graph $(V(K_n), E_{j, t})$, where $E_{j, t}$ is the set of edges taken by the $j$-th player up to his $t$-th turn. We also let $\Gamma_t$ be the graph $(V(K_n), E_{1,t} \cup E_{2,t})$. We will usually abuse notation by omitting $t$ when the turn is understood.

Using this notation, we have that $\Delta(H_j) \leq k$ for $j \in \{1, 2\}$ and hence $\Delta(\Gamma) \leq 2k$ holds before any player loses. In particular, the game cannot be a draw for $n \geq 2k + 2$. 
\\\\
\textit{Proof of Theorem 1.} Let $n, k$ be positive integers such that $n\geq 200k$ and consider the $(k + 1)$-star avoidance game on $K_n$. Suppose for contradiction that PII doesn't win the game. The following claim is key to the proof:
\\\\
\textbf{Claim.} \textit{PII can ensure that eventually one or two vertices in $H_2$ have degree $k - 1$ while the rest have degree $k$, and additionally the vertices of degree $k - 1$ span no edges in $\Gamma$.}\\

We first show how the Claim implies the desired result. Note that PI certainly loses on his $\left(\lfloor \frac{nk}{2} \rfloor + 1 \right)$-st move. Let $\Phi$ denote the strategy provided by the Claim. If $nk$ is odd, PII follows $\Phi$ and hence wins. So suppose $nk$ is even. Then PII follows $\Phi$ for the first $\frac{nk}{2} - 2$ rounds. After PI's $\left(\frac{nk}{2} - 1\right)$-st move, let $e_1, e_2$ denote PII's next move according to $\Phi$ and the pair of vertices that would have degree $k - 1$ in $H_2$ if PII claimed $e_1$, respectively. Note that PI's next move is fixed at this moment, so if it is among $\{e_1, e_2\}$, PII simply wins by claiming it. Otherwise, PII claims $e_1$ and $e_2$ in that order and hence wins.
\\\\
\textit{Proof of Claim.} PII's strategy is divided into $k$ stages. For all $j$, at the beginning of the $j$-th stage, at most $2$ vertices in $H_2$ will have degree $j$ and the rest will have degree $j - 1$. We will show by induction on $j$ that PII will be able to maintain this property. So fix some $j \in [k]$ and throughout the $j$-th stage, let $S = \{v \in V(H_2) \mid d_{H_2}(v) = j - 1\}$. PII follows the winning strategy for PCG($G$), where $G = \Gamma[S]$. By Theorem 5, he can indeed do so since $v(G) \geq n - 2 \geq 200k - 2$ and $\Delta(G) \leq k + j - 1 \leq 2k - 1$ hold at the beginning. Consequently, we are done in the case when $j = k$. Otherwise, if $j < k$, then as long as $S$ is non-empty, PII chooses a vertex in $S$ and a vertex of degree $j$ in $H_2$. He can do so since $\Delta(\Gamma) \leq 2k$ and there are at least $n - 2 \geq 200k - 2$ vertices of degree $j$ in $H_2$. This results in at most $2$ vertices of $H_2$ having degree $j + 1$ and the rest having degree $j$, as desired. $\blacksquare$

\bigskip 

\section{Concluding remarks and open problems}

As a corollary of our result, there exists a function $h : \mathbb{N} \to \mathbb{N}$ with the property that $h(k)$ is the least positive integer $n_0$ such that the $(k + 1)$-star avoidance game on $K_n$ is a PII win for all $n \geq n_0$. Using this notation, Theorem 1 may be rephrased as the assertion that $h(k) \leq 200k$ holds for all positive integers $k$. It is clear that this bound is not optimal, in particular it is possible to verify with the help of a computer that $h(1) = 3$ and $h(2) = 5$. We are also aware that certain modifications of the presented approach may lead to replacing $200$ by a smaller constant. As has already been remarked in Section 2, the game is not a draw for $n \geq 2k + 2$. Hence, we find it natural to ask the following question:
\\\\
\textbf{Question 6.} \textit{Is it true that $h(k) \leq 2k + 2$ for all positive integers $k$?}\\

We have little intuition as to the correct answer to this question. In particular, we doubt that our approach can be modified so as to settle this question. The reason for this is that Question 1 having an affirmative answer would probably have to do with the game having a high degree of symmetry rather than it being very 'sparse', in the sense that ex$(n, S_{k+1})$ is much smaller than $\binom{n}{2}$ for large $n$. It would also be interesting to find a complete characterisation of the outcomes of the game:
\\\\
\textbf{Question 7.} \textit{For each pair of positive integers $(k, n)$, is the $(k + 1)$-star avoidance game on $K_n$ a PI win, a PII win or a draw?}\\

A related problem would be to consider sim-like games with certain classes of sparse graphs instead of $S_{k+1}$, for example trees, and in particular paths. It would be interesting to pursue results of a similar kind for such games. Finally, no example of a sim-like game that is a PI win is known. The following question of Johnson, Leader and Walters (see \hyperlink{transitive}{[5]}) remains open:
\\\\
\textbf{Question 8.} \textit{Does there exist a sim-like game that is a PI win?}

\section*{Acknowledgments}

I am grateful to Prof. Imre Leader for his support and useful discussions, as well as suggestions that greatly improved the presentation of the paper.
\\\\
I was provided financial support by Trinity College, Cambridge, for which I am also grateful.

\end{document}